\documentclass[compress,final,3p,times,11pt]{elsarticle}
\usepackage{amsfonts}
\usepackage{amssymb}
\usepackage{amsmath}
\usepackage{color,epstopdf}
\usepackage{mathrsfs}
\usepackage{graphicx,color,epstopdf}
\usepackage{float}
\usepackage{caption}
\usepackage{subcaption}
\usepackage{float}
\usepackage{subcaption}

\newcommand{\al}{X}
\journal{Journal of Statistical Mechanics: Theory and Experiment}
\begin{document}

\title{Stochastic Bifurcation in Single-Species Model Induced by $\alpha$-Stable L\'evy Noise \tnoteref{mytitlenote}}
\author[mymainaddress,mysecondaryaddress]{Almaz Tesfay}
\ead{amutesfay@hust.edu.cn}
\author[mymainaddress,mysecondaryaddress]{Daniel Tesfay}
\ead{dannytesfay@hust.edu.cn}
\author[mymainaddress]{Shenglan Yuan\corref{mycorrespondingauthor}}
\ead{shenglanyuan@hust.edu.cn}
\author[mythirdaddress]{James Brannan}
\ead{jrbrn@clemson.edu}
\cortext[mycorrespondingauthor]{Corresponding author}
\author[myfourthaddress]{Jinqiao Duan}
\ead{duan@iit.edu}
\address[mymainaddress]{School of Mathematics and Statistics \& Center for Mathematical Sciences,  Huazhong University of Science and Technology, Wuhan 430074, China}
\address[mysecondaryaddress]{Department of Mathematics, Mekelle University, P.O.Box 231, Mekelle, Ethiopia}
\address[mythirdaddress]{Department of Mathematical Sciences, Clemson University,Clemson, South Carolina 29634, USA}
\address[myfourthaddress]{Department of Applied Mathematics, Illinois Institute of Technology, Chicago, IL 60616, USA}
\begin{abstract}
Bifurcation analysis has many applications in different scientific fields, such as electronics, biology, ecology, and economics. In population biology, deterministic methods of bifurcation are commonly used. In contrast, stochastic bifurcation techniques are infrequently employed. Here we establish stochastic P-bifurcation behavior of (i) a growth model with state dependent birth rate and constant death rate, and (ii) a logistic growth model with state dependent carrying capacity, both of which are driven by multiplicative symmetric stable L\'evy noise. Transcritical bifurcation occurs in the deterministic counterpart of the first model, while saddle-node bifurcation takes place in the logistic growth model. We focus on the impact of the variations of the growth rate, the per capita daily adult mortality rate, the stability index, and the noise intensity on the stationary probability density functions of the associated non-local Fokker-Planck equation. Implications of these bifurcations in population dynamics are discussed. In the first model the bifurcation parameter is the ratio of the population birth rate to the population death rate. In the second model the bifurcation parameter corresponds to the sensitivity of carrying capacity to change in the size of the population near equilibrium. In each case we show that as the value of the bifurcation parameter increases, the shape of the steady-state probability density function changes and that both stochastic models exhibit stochastic P-bifurcation. The unimodal density functions become more peaked around deterministic equilibrium points as the stability index increases. While, an increase in any one of the other parameter has an effect on the stationary probability density function. That means the geometry of the density function changes from unimodal to flat, and its peak appears in the middle of the domain, which means a transition occurs.

\end{abstract}
\begin{keyword}
 stochastic bifurcation; Brownian motion; L\'evy motion; single-species model; stochastic differential equation;.
\MSC[2020] 39A50, 45K05, 65N22.
\end{keyword}
\maketitle
\section{Introduction}
Population biology studies the dynamics and changes in the size of populations. The dynamics include the growth and decline of populations and the interaction of species. In modeling the dynamics, different mathematical population models have been investigated,   many in the form of ordinary differential equations.

We consider a one-dimensional deterministic model of the population for a certain biological species
 \cite{1,2,4,3}
\begin{equation}\label{equat2}
\frac{d\al}{dt}=\al\,f(\al).   \tag {1.1}
\end{equation}
The function $f(\al)$ is the per capita growth rate of the population which depends on on  the population size $\al$. According to \cite{5}, $f(\al)$ is also called net per capita growth rate. Different choices of $f(\al)$ in Eq. (\ref{equat2}) give us different population growth equations.

The choice  $f(\al)=(-p+qe^{-s\,\al})$, gives the model \cite{6,7}
\begin{equation}\label{deter1}
\frac{d\al}{dt}=\al\,(-p+h(\al)),\tag {1.2}
\end{equation}
where $h(\al)=q\,e^{-s\,\al}$ denotes the population dependent birth rate, $s >0$  represents the per capita daily adult mortality rate, $q > 0$ is the birth rate when the population is small, $p>0$ is the death rate with $q > p$. The population size $X$ at time $t$ in a given place, is always positive ( cf. Section \ref{Sec3}).

If we set $f(\al)=s\,(1-\frac{\al}{M})$, then Eq. (\ref{equat2}) becomes the Verhulst model (logistic growth equation) \cite{5,DK,DK19,aad,8}
 \begin{equation}\label{Det21}
\frac{d\al}{dt}=s\,\al\,\left(1-\frac{\al}{M}\right), \quad \al_0=x_{0},\tag {1.3}
 \end{equation}
where $\al$ is the population size at time $t$ in a given area or place, $s> 0$ is the population growth rate when the population is small, and carrying capacity $M >0$ is usually determined by the available sustaining resource.  The carrying capacity is often assumed to be a constant \cite{9,10}. In this paper we assume a state-dependent carrying capacity $M(X)$( cf. Section \ref{Sec4}).

Models (\ref{deter1}) and (\ref{Det21}) efficiently describe the biological growth depending on environmental factors such as food, flood, and weather.

The deterministic version of the logistic model, a nonlinear differential equation, is used to model biological populations in which the per capita growth rate $f(X)$ decreases linearly as the population size $X$ approaches a maximum characteristic size $M$, and
in fact turns negative for $X > M$ . The logistic model has long been applied in ecological modeling because of its ability to simply capture dominant features, early term exponential growth and late term boundedness on growth, while requiring only
two parameters \cite{11}.

Stochastic differential equations (SDEs) are important tools for modeling and analysis of biological population growth systems \cite{12}. Numerous environmental factors affect  population growth: temperature, earthquake, floods and so on, need to be considered. By adding stochasticity to the models, it is possible to study the effects of such perturbations on population growth that deterministic models can not account for.

Next, let us see how stochastic fluctuations of population growth are added to population models. Frequently, environmental fluctuations are modeled by adding noise to the per capita growth rate $f(\al)$;  see \cite{aad1,13,aad2}. In other words, the stochastic version of models  (\ref{deter1}) and (\ref{Det21}) may be constructed by adding \emph{Gaussian noise} or \emph{non-Gaussian noise} or \emph{both noises} to $f(X)$. We add noise to the specific rate which emerge representing environmental fluctuations which chiefly affect growth or reproduction of individuals belonging to the population.

\begin{align}\label{eq1-SDE}
f(X)\rightarrow f(X)+\lambda\,\dot{B}_t+ \epsilon\,\dot{L}_t^{\alpha}.\tag {1.4}
\end{align}

The influence of the environment on the subsystem is then accounted for in a population sense by treating coefficients and/or input to Equations (\ref{deter1}) and (\ref{Det21}) as random variables or processes whose statistical properties are supposed to be known. The solution of the equations will be a random process, and the problem consists of finding their statistical properties as well as the statistical properties of certain functional of the solution.

In Eq. (\ref{eq1-SDE}), the function $f(\al)$ represents for a mean fractional growth rate, the stochastic differential equation of Eq. (\ref{equat2}) has the following form:
\begin{equation}\label{eq1-SDE2}
\frac{d\al}{\al}=f(\al)dt + \lambda dB_t + \epsilon dL_t^{\alpha},\tag {1.5}
\end{equation}
or equivalently,
\begin{equation*}
d\al=\al\,f(\al)dt + \lambda\,\al dB_t + \epsilon\,\al dL_t^{\alpha}.
\end{equation*}
Here, we have used standard upper case notation for the process $\al=\{X_t, t\geq 0\}$.\par
In this study, we assume the mean fractional growth rate to be of the form $f(\al)=(-p+h(\al))$ and $f(\al)=s\, \left(1-\frac{\al}{M}\right)$. This gives the standard stochastic population growth equations driven by combined Gaussian and non-Gaussian noise. The SDE version of single-species becomes
\begin{equation}\label{deter1-SDE}
d\al=\al\,(-p+h(\al))dt+\lambda\, \al\, dB_t+\epsilon\, \al \,dL_t^\alpha,\qquad \al_0=x_{0},\tag {1.6}
\end{equation}
and the logistic growth SDE version is defined by
 \begin{equation}\label{Det21-SDE}
d\al=s\,\al\, \left(1-\frac{\al}{M}\right)dt +\lambda \, \al \,dB_t+\epsilon\, \al\, dL_t^\alpha , \qquad \al_0=x_{0}.\tag {1.7}
 \end{equation}
 In both Eqs. (\ref{deter1-SDE}) and (\ref{Det21-SDE})  $B_t$ is the one dimensional Brownian motion, the parameter $\lambda\in (0,1)$ is the intensity of Gaussian noise,  $\epsilon\in (0,1)$ represents non-Gaussianity index, $L_{t}^{\alpha}$ is an $\alpha$-stable L\'evy motion with stability index $\alpha\in (0,2)$

Recently, Sun et. al. \cite{6} studied Eq. (\ref{deter1-SDE}) in the Gaussian noise case only. The Authors  proved the existence and uniqueness of the solution, and they also discussed the recurrence and the existence of stationary distribution. In our research we consider Eq. (\ref{deter1-SDE}) in the non-Gaussian noise case and study about the effect of parameters in the model.

Our previous papers \cite{9,10} assumed a constant carrying capacity in Eq. (\ref{Det21-SDE}) and investigated the extinction probability of the population. But here we assume the state dependence maximum population size $M(X)$ and we focus on the stochastic bifurcation of the model.

In this work, we mainly focus on deterministic bifurcation and use stationary probability density functions to describe stochastic P-bifurcation of each for the two population growth models given in (\ref{deter1-SDE}) and (\ref{Det21-SDE}) one by one. More details are given in Sections  \ref{Sec3} and \ref{Sec4}.

Many scholars are interested in stochastic bifurcations \cite{14,15,17,18,16}. There are few rigorous general theorems and criteria to detect stochastic bifurcations, which are often only verified by computer simulations or for some particular models \cite{19}. Stochastic bifurcation problems are more difficult than deterministic bifurcation problems. Stochastic bifurcations can be described by a qualitative transition of the stationary probability distribution, for instance, a change from unimodal to a flatened shape, or from unimodal to a bimodal distribution. According to \cite{20,18}, stochastic bifurcations are divided into two kinds:  phenomenological P-bifurcation (based on the change of shape of the stationary probability density function), and dynamical D-bifurcation based on the change of sign of the largest Lyapunov exponential. D-bifurcation is a dynamic concept, similar to deterministic bifurcations, while P-bifurcation is a static concept. Unfortunately, these two definitions do not agree well, which implies a need for exploring a new definition of stochastic bifurcation.

Stochastic bifurcations have been investigated in a broad scope of nonlinear systems in engineering and physical science.  Several authors have studied stochastic bifurcation for stochastic dynamical systems with Brownian motion  by detecting  the effects of varying parameters in stationary probability densities as solutions of Fokker-Planck equations. \cite{21} obtained some properties of the random Conley index and then presented a sufficient condition for the existence of abstract bifurcation points for continuous-time and discrete-time random dynamical systems.  \cite{22} measured the critical bifurcation parameter for a first-order phase transition in non-equilibrium systems in the thermodynamic limit and evaluated the mechanism of phase transition.

Based on stochastic bifurcation for dynamical system driven by non-Gaussian $\alpha$-stable L\'evy noise \cite{YZD}, in this study, we are motivated to explore the dynamical behaviors of  population growth models with L\'evy stable noise perturbation. This helps to demonstrate the interactions between the population models and their complex surrounding. Here, we focus on the stochastic bifurcation  of the single-species  model driven $\alpha$-stable L\'evy jump noise affected by parameters and discuss biological interpretations of our findings.

In the present study, we investigate P-bifurcation for a type of growth and the Verhulst model with state-dependent carrying capacity in the presence of $\alpha$-state L\'evy noise. We think this is a better tool to demonstrate the interactions between the population system and its complex surrounding. Especially, we study stochastic bifurcation induced by the deterministic parameters and L\'evy stable noise, and present the biological interpretations. The existence of the equilibrium solutions of the two deterministic models (\ref{deter1}) and (\ref{Det21}) is analyzed. We determine the deterministic bifurcation of Eqs. (\ref{deter1}) and (\ref{Det21}) and stochastic bifurcation  of the systems (\ref{deter1-SDE}) and (\ref{Det21-SDE}).

This paper is structured as follows. We begin with the basic concepts about Brownian motion, L\'evy $\alpha$-stable motion, and non-local Fokker-Plank equation in Section \ref{Sec2}. We review the deterministic counterpart and analyze the stochastic bifurcation of a single-species population model driven by multiplicative non-Gaussian noise in Section \ref{Sec3}. In Section \ref{Sec4}, we study the stochastic logistic growth model with state variable carrying capacity under multiplicative non-Gaussian L\'evy noise. Moreover, we investigate the deterministic and stochastic bifurcations that underlie qualitative change of the stationary probability density. Section \ref{Sec5} discusses our numerical results in two models above. Finally, we finish this research with a conclusion and discussion in Section \ref{Sec6}.
\section{Preliminaries}\label{Sec2}
In this section, we review some basic facts about Brownian motion, L\'evy motion and  non-local Fokker Planck equation.
\subsection{\textbf{{Brownian motion}}}\label{Sec22}
The one-dimensional standard Brownian motion $B_t$ is stochastic process, defined on complete probability space $(\Omega, \mathfrak{F},\mathfrak{F}_t,\mathbb{P})$; see \cite{24}. A Gaussian zero-mean process $B_t$ having the following properties: (i) independent and stationary increments; (ii) the variance of $B_t-B_s\, (t\geq s \geq 0)$ is $t-s$ and $B_0=0$ almost surely (a.s); (iii) continuous sample paths almost surely  and its paths are nowhere differentiable. It plays an important role in describing many physical phenomena that exhibit random behavior.
\subsection{\textbf{The $\alpha$-Stable L\'evy motion}}\label{Sec23}
L\'evy motions $L_t$ are a class of non-Gaussian stochastic processes that also have independent and stationary increments. A stable distribution $S_{\alpha}(\theta,\beta,\gamma)$ is the distribution for a stable random variable, where the stability index $\alpha \in (0,2)$, the skewness $\beta \in (0,\infty)$, the shift $\gamma \in (-\infty,\infty)$, and scale index $\theta \geq0$. An $\alpha$-stable L\'evy motion $L_{t}^{\alpha}$ \cite{25,26,27,28} is a non-Gaussian stochastic process satisfy (i) $L_{0}^{\alpha}=0$, almost surely;
(ii) $L_{t}^{\alpha}$ has independent increments: the random variables $L_{t_{i+1}}^{\alpha}-L_{t_{i}}^{\alpha}$ are independent for $0\leq t_{1}<t_{2}<t_{3}<...<t_{i-1}<t_{i}<\infty,$ and for each $i=1,2,...$;
(iii) $L_{t}^{\alpha}$ has stationary increments: $L_{t}^{\alpha}$-$L_{s}^{\alpha}$ and $L_{t-s}^{\alpha}$ have the same distribution $S_{\alpha}((t-s)^{1/{\alpha}},0,0)$;
(iv) $L_{t}^{\alpha}$ has stochastically continuous sample paths, i.e., for $0 \leq s\leq t$ and $\delta>0$, the probability $\mathbb{P}(|L_{t}^{\alpha}-L_{s}^{\alpha}|>\delta) \rightarrow 0$ as $t\rightarrow s$.

An $\alpha$-stable L\'evy motion $L_t^{\alpha}$ having values in $\mathbb{R} = (-\infty, \infty)$ is defined by a drift coefficient $\hat{b}\in \mathbb{R}$, $\hat{Q}\geq0$ and a Borel measure $\nu_{\alpha}$ defined on ${\mathbb{R}}\setminus{\{0\}}$. The triplet $(\hat{b},\hat{Q},\nu_{\alpha})$ is the so-called generating triplet of L\'evy motion $L_t^{\alpha}$.  A L\'evy motion can be expressed as linear combination of time $t$, a Brownian motion $B_t$ and a pure jumping
process \cite{30,29}, i.e., $L_t^{\alpha}$ can be written as
\begin{equation}\label{Ito-Dec}
 L_t^{\alpha}=\hat{b}t+ B_{\hat{Q}}(t)+\int_{|y|<1}y\tilde{N}(t,dy)+\int_{|y|\geq 1}yN(t,dy),\tag {2.1}
\end{equation}
where $N(t,dy)$ is the independent  Poisson random measure on $\mathbb{R}^+\times{{\mathbb{R}}\setminus{\{0\}}}$, $\tilde{N}(t,dy)=N(t,dy)-\nu_{\alpha}(dy)dt$ is the compensated Poisson random measure, $\nu_{\alpha}(S)=E(N(1,S))$ is the jump measure, and $B_t$ is the an independent standard 1-dimensional Brownian motion.

The L\'evy-Khinchin formula for $\alpha$-stable L\'evy motion has a specific form for its characteristic function. For $0\leq t<\infty$, $u \in \mathbb{R}$,
\begin{equation*}
\mathbb{E}[e^{(iu{L_t^{\alpha}})}] = e^{(t\psi(u))},
\end{equation*}
where
$$\psi(u)=iu\hat{b}-\frac{\hat{Q}}{2}u^2 + \int_{{\mathbb{R}}\setminus{\{0\}}}(e^{iuz}-1-iuzI{_{|z|<1}})\nu_{\alpha}(dz).$$
L\'evy measure $\nu_{\alpha}$ is given by
 $\nu_{\alpha}(du)=c({\alpha})\frac{1}{|u|^{1+\alpha}}du$
with $ c(\alpha) = \alpha\frac{\Gamma(\frac{1+\alpha}{2})}{{2^{1-\alpha}\pi^{\frac{1}{2}}}\Gamma{(1-\frac{\alpha}{2})}}$, and $\Gamma$ is the Gamma function. The function of the L\'evy measure is to describe the expected number of jumps of a certain size at a time interval one. Usually, the parameter $\alpha$ is called the index of stability with the value $0 < \alpha < 2$.
In the case of a one-dimensional $\alpha $-stable L\'evy motion, the drift vector $\hat{b}=0$ and the diffusion $\hat{Q}=0$.
Here, we focus on jump process with a specific size in generating triplet $(0,0,\nu_{\alpha})$ for the random variable $S_{\alpha}$ which can be defined by $\Delta L_{t}^{\alpha}=L_{t}^{\alpha}-L_{t^-}^{\alpha}<\infty, t \geq 0,$ where $L_{t-}^{\alpha}$ is the left limit of the L\'evy motion in $ R=(-\infty, \infty)$ at any time $t$.

\subsection{Fokker-Planck equation }
The Fokker-Planck equation (FPE) is a deterministic equation that describes how the conditional probability density function of the solution process evolves as time progresses.

Consider a scalar stochastic differential equation driven by both Gaussian noise and non-Gaussian noise:
\begin{equation}\label{Gen-SDE}
d\al=f_1(\al)\,dt + f_2(\al)\,dB_t + f_3(\al)\,dL_t^{\alpha},\tag {2.2}
\end{equation}
where $f_1(\al)=\al\,f
(\al)$ is a given drift term, $f_2(\al)=\lambda\,\al$ is the diffusion coefficient and $f_3(\al)=\epsilon\,\al$ is L\'evy noise intensity, $B_t$ is the standard Brownian motion, $\lambda$ is the intensity of Gaussian noise, $\epsilon$ represents non-Gaussianity index, and $L_{t}^{\alpha}$ is an $\alpha$-stable L\'evy motion with stability index $\alpha\in (0,2)$.

Define $z=f_3(\al)y$. The generator $\mathcal{A}$ of the stochastic differential equation in (\ref{Gen-SDE}) \cite{32,31,33} is given by
\begin{align}\label{Generator}
\mathcal{A}\,\varphi(\al)= f_1(\al)\,\varphi'(\al)+\frac{1}{2}{f_2^2(\al)}\,\varphi''(\al)
  +|f_3(\al)|^{\alpha}\,\int_{\mathbb{R}\setminus{\{0\}}}\left[\varphi(\al +z)-\varphi(\al)\right]\nu_\alpha(dz). \tag {2.3}
 \end{align}

\subsubsection*\textbf{{Non-local Fokker-Planck equation:}}

The $P(\al,t)$ is probability density function for the solution of Eq. (\ref{Gen-SDE}) with $\lambda=0$ subject to the initial condition $\al_0=x_0$, the non-local Fokker-Planck equation for $P(\al,t)$  is \cite{35,34,36}
\begin{equation*}
 P_t(\al,t)=\mathcal{A}^{*}P(\al,t), \qquad \al \in D,
 \end{equation*}
\begin{equation*}
P(\al,0)=\delta(\al-x_0), \qquad x_0\in D,
\end{equation*}
where $\delta $ is the  Dirac function. The ad-joint operator $\mathcal{A}^*$ of $\mathcal{A}$ in the, Hilbert space $L^2(R)$ is obtained by solving
\begin{equation*}
\int_{\mathbb{R}\setminus{\{0\}}}\mathcal{A}\varphi (\al)V(\al)d\al=\int_{\mathbb{R}\setminus{\{0\}}}\varphi (\al)\mathcal{A}^*V(\al)d\al,
\end{equation*}
for $\varphi $, $V$ in the domain of definition for the operator $\mathcal{A}$ and $\mathcal{A}^*$, respectively. We find that
 \begin{equation*}
\mathcal{A}^*V(\al)=\int_{\mathbb{R}\setminus{\{0\}}}[|f_3(\al)+z)|^{\alpha}V(\al+z)-|f_3(\al)|^{\alpha}V(\al)]\nu_\alpha(dz).
\end{equation*}
Therefore, according to \cite{33}, $P(\al,t)$ satisfies
\begin{align}\label{FPE-Levy}
\frac{\partial P}{\partial t}&=-\frac{\partial}{\partial \al}\,(f_1(\al)P(\al,t))+\frac{1}{2}\frac{\partial ^2}{\partial \al^2}(f_2^2(\al)P(\al,t))\nonumber\\
&+|\epsilon\,\al|^{\alpha}\,\int_{\mathbb{R}\setminus{\{0\}}}[|f_3(\al+z)|^{\alpha}P(\al+z,t)-|f_3(\al)|^{\alpha}P(\al,t)]\nu_\alpha (dz).\tag {2.8}
\end{align}

Since it is very difficult to find the analytical solutions of the non-Fokker-Planck equation in Eq. (\ref{FPE-Levy}). We find numerical approximations using the simulation method developed in \cite{35}.
\section{ Bifurcation analysis of a growth model}\label{Sec3}
\subsection{\textbf{Deterministic transcritical bifurcation}}\label{mod1}
This section focuses on the dynamics of the generalized form of the deterministic growth model. Recall Eq. (\ref{deter1})
\begin{equation}\label{deter}
\dot{\al}=\al\,(-p+h(\al)),\tag {3.1}
\end{equation}
where $h(\al)=q\,e^{-s\,\al}$ represents the population dependent birth rate, $q > 0$ is the birth rate when the population is small, $s >0$  denotes the per capita daily adult mortality rate, $p>0$ is the death rate with $q > p$ satisfying the following conditions:
\begin{description}
   \item[i.] the function $h(\al)$ is positive, i.e., $h(\al) > 0$;
  \item[ii.] the function $h(\al)$ is continuously differentiable with $h'(\al) < 0$, and $p$ is between $ h(0^+)$ and  $h(\infty)$.
\end{description}
The equilibrium solutions of Eq. (\ref{deter}) are $\al_1=0$ and $\al_2= \frac{\ln q-\ln p}{s}=\frac{1}{s}\ln\frac{q}{p}=\frac{1}{s}\ln\mu,$ where $\mu=\frac{q}{p}$.
\begin{figure}[ht!]
    \centering
     \begin{subfigure}[b]{0.45\textwidth}
         \centering
         \includegraphics[width=\textwidth]{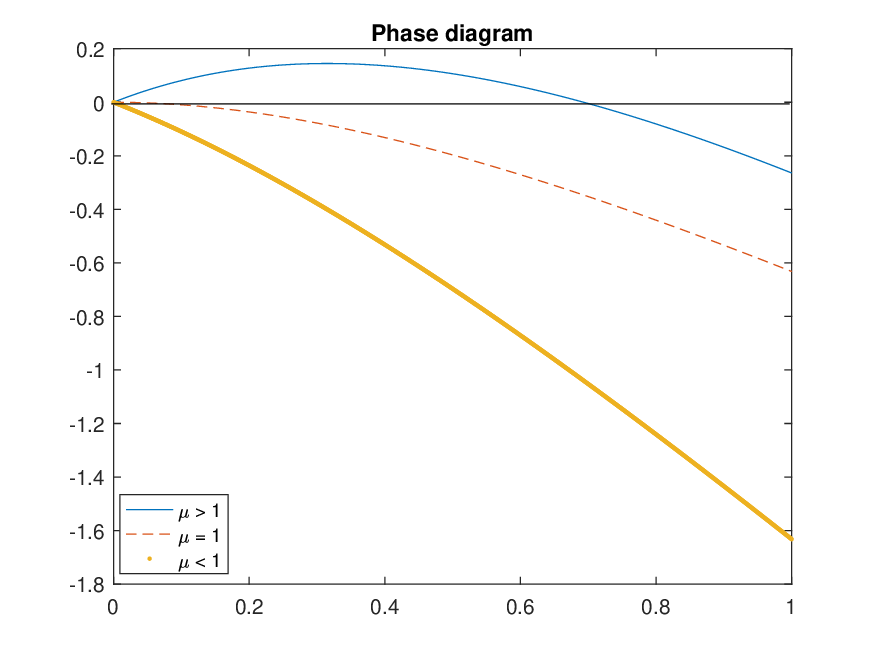}
         \caption{$\frac{d\al}{dt}$ verse $\mu$.}
          \end{subfigure}
     \begin{subfigure}[b]{0.45\textwidth}
         \centering
         \includegraphics[width=\textwidth]{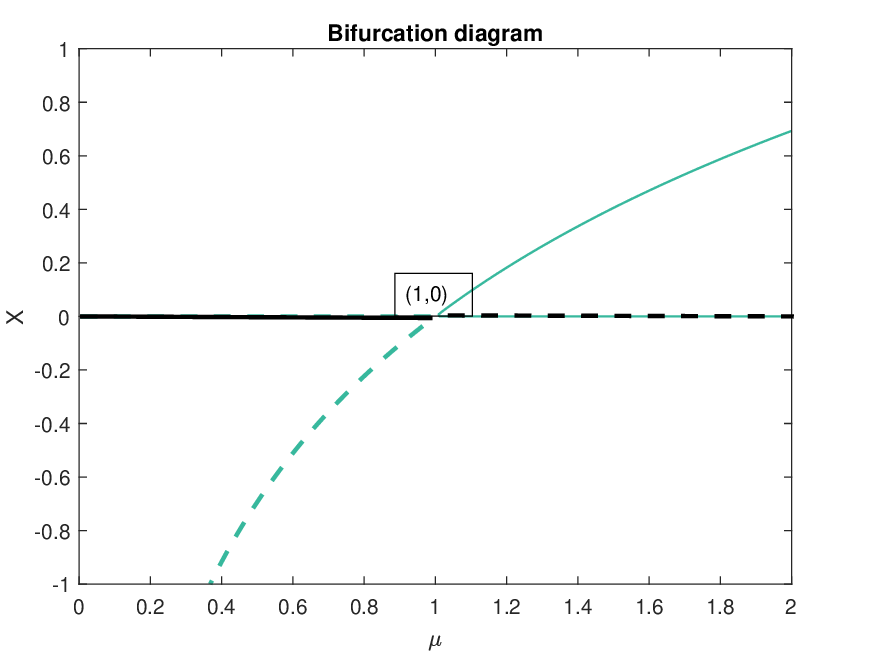}
         \caption{$\mu=\frac{q}{p}$.}
         \end{subfigure}
     \begin{subfigure}[b]{0.55\textwidth}
         \centering
         \includegraphics[width=\textwidth]{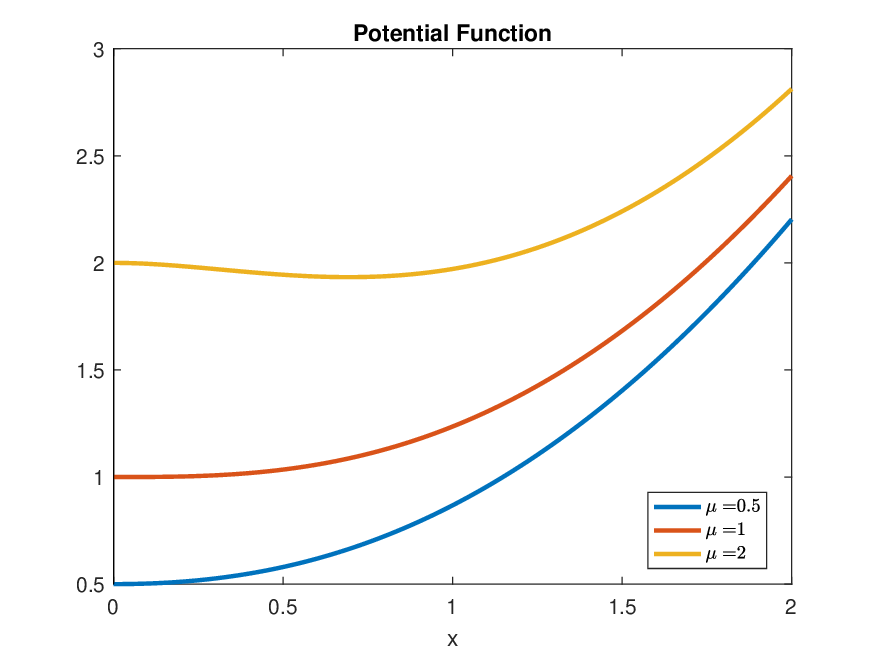}
         \caption{Potentian function $U(\al)$ verse $\mu$.}
         \end{subfigure}

        \caption{ (Color online) (a) Phaseline diagram of Eq. (\ref{deter}) when different values of $\mu$ with fixed $s=1$. (b) The bifurcation diagram of  Eq. (\ref{deter}) for $s=1$. The broken line is unstable state and the solid line is stable state. (c) The potential function model (\ref{deter}). The potential has its minimum at $\al_2$, and  maximum at $\al_1$.}
 \label{Ph-M1}
\end{figure}

If there is a qualitative change in the dynamics upon a small change in a parameter, then we say a bifurcation occurs. The simplest bifurcations correspond to changes in equilibria, namely, in their number and stability type. For model (\ref{deter}), there are two equilibria when $p\neq q$, and one when $p=q$. Thus there is a bifurcation, a change in the number of equilibria, on the line $p=q$ in the parameter space; this is the bifurcation set. The existence of the equilibria depends only on one combination
of the two parameters, $\mu=q/p$; consequently this bifurcation is governed by a single effective parameter. We can conveniently collect the information about the equilibria in a bifurcation diagram that shows the two functions $\al_1$ and $\al_2$ of the single parameter $\mu$.

Denote the right side of Eq. (\ref{deter}) by $F(\al)$. Then $F'(\al)=-p+(1-s\al)\,h(\al)$. Hence $F'(\al_1)=q-p$ and $F'(\al_2)=-p\ln \mu$. For $\mu<1$, $\al_1=0$ is stable and $\al_2=\frac{1}{s}\ln \mu$ is unstable. These two equilibria coalesce at $\mu=1$. For $\mu>1$,  $\al_1=0$ is unstable and $\al_2=\frac{1}{r}\ln \mu$ is stable. Thus, an exchange of stability occurs at $\mu=1$. The bifurcation occurs at the point $(\mu,\al)=(1,0)$ where the two equilibria collide. According to the book \cite[ pg=26]{37}, this type of bifurcation is called a transcritical bifurcation.

In the deterministic model (\ref{deter}), two curves of equilibria intersect at the point $(1,0)$ in the bifurcation diagram. Both curves exist on the either side of $\mu=1$. However, the stability of the equilibrium solution along a given curve changes on passing through $\mu=1$.

By substituting the Taylor expansion of $e^{-s\,\al}$ into Eq. (\ref{deter}) and ignoring the higher orders with small population size, Eq. (\ref{deter}) has normal form for a transcritical bifurcation \cite{38} which is given by $$\dot{\al}=F(\al,\lambda)=\lambda\, \al-\al^2,$$ where $\lambda=\frac{q-p}{\sqrt{s\,q}}$. The bifurcation diagram (the $\lambda-\al$ plane) represents the qualitative behavior of the single-species model (\ref{deter}). In this Figure, two curves of equilibria ($\al_1=0$ and $\al_2=\lambda$) intersect at the point $(0,0)$. That shows that it has the bifurcation point (0,0) where a stable and an unstable equilibrium collide and stability changes. This type of bifurcation is called a ''transcritical" bifurcation (exchange of stabilities) \cite[pg=360]{39}.

Model (\ref{deter}) can be further expressed as its potential function $U$, i.e.,
 \begin{equation*}
\frac{d{\al}}{dt}=-\frac{dU}{d\al},
\end{equation*}
where
\begin{align*}
U(\al)= &-\int\al\,(-p+q\,e^{-s\,\al})d\al\nonumber\\
  = &\frac{p}{2}\al^2+\frac{q}{s}\,e^{-s\,\al}[\al+\frac{1}{s}]\nonumber\\
  = &\mu\left[\frac{p^2}{2q}\al^2+\frac{p}{s}e^{-s\,\al}(\al+\frac{1}{s})\right],  \quad for\, \,  \mu=\frac{q}{p}.
\end{align*}
The plot in Fig. \ref{Ph-M1}c shows the potential function $U$ of the deterministic equation (\ref{deter}) with different values of $\mu$.
\subsection{\textbf{Stochastic bifurcation for a growth model} }\label{sec32}
\begin{figure}[htb!]
\begin{center}
  \includegraphics[width=0.5\linewidth]{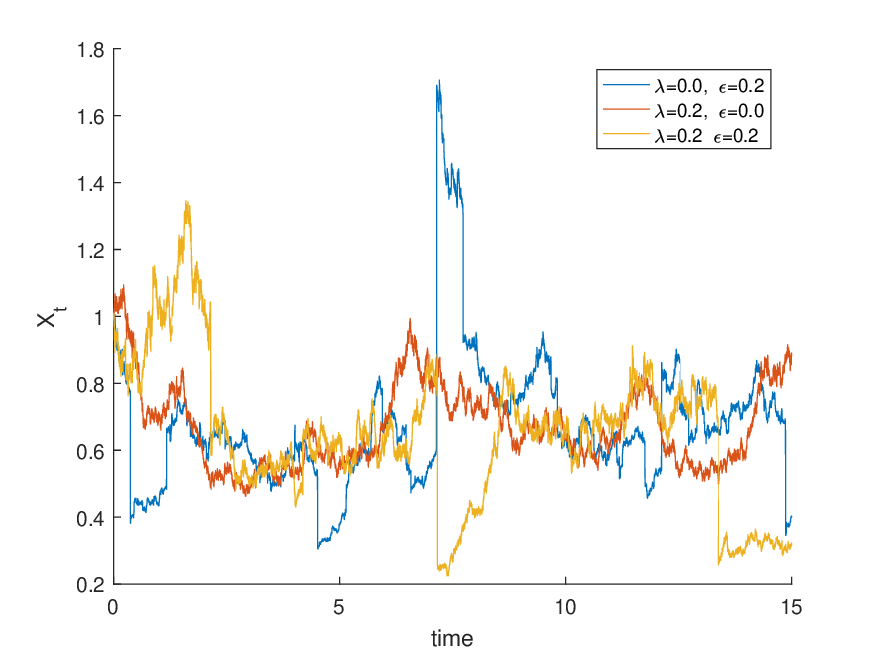}
 \caption{The numerical solution of SDE  (\ref{s1}) for $s=1,\,\, q=2,\,\,p=1$, and $x_0=1$.}
 \label{sample-M1}
 \end{center}
\end{figure}

In this subsection, we consider the stochastic population model:
\begin{equation}\label{s1}
d\al=\al\,(-p + q\,e^{-s\,\al})dt +\lambda\, \al\,dB_t+\epsilon \,\al \,dL_{t}^{\alpha},\quad \al_0=x_0,\tag {3.2}
\end{equation}
where $B_t$ is the standard Brownian motion, $\lambda$ is the intensity of Gaussian noise and $\epsilon$ represents the intensity of non-Gaussian noise, $L_{t}^{\alpha}$ is an $\alpha$-stable L\'evy motion, and $\alpha \in (0,2)$ is stability index. Here the noises are multiplicative because the diffusion and the intensity coefficient in (\ref{s1}) depend on $\al$. Figure \ref{sample-M1} shows a sample path of the SDE in (\ref{s1}) when $\epsilon=0$.
Next, by assuming $\al$ in (\ref{s1}) has a conditional probability density $p(\al,t):=p(\al,t|x_{0},0)$, we derive Fokker-Plank equation satisfied by $p(\al,t)$.

From Eq. (\ref{Generator}) the generator $\mathcal{A}$ for the model (\ref{s1}), or for the solution process $\al$, is
\begin{equation*}
\mathcal{A}\,\varphi(\al)=\al\,(-p+q\,e^{-s\al})\,\varphi'(\al)+\frac{1}{2}\lambda^{2}\al^{2}\varphi''(\al)+|\epsilon\,\al|^{\alpha}\,\int_{\mathbb{R}\setminus{\{0\}}}\left[\varphi(\al +z)-\varphi(\al)\right]\nu_\alpha(dz),
\end{equation*}
where $\varphi$ belongs to the domain of definition for the operator $\mathcal{A}$. The Fokker-Planck equation for the SDE (\ref{s1}) in terms of the probability density function for the solution process $\al$ with initial condition $\al_0=x_{0}$, is
\begin{align}\label{p1}
\partial_{t}P(\al,t)&=-\partial_{\al}(\al\,(-p+q\,e^{-s\al})\,P(\al,t))+\frac{\lambda^{2}}{2}\partial_{\al\,\al}(\al^{2}P(\al,t))\nonumber\\
&+|\epsilon\,\al|^{\alpha}\,\int_{\mathbb{R}\setminus{\{0\}}}[|f_3(\al+z)|^{\alpha}P(\al+z,t)-|f_3(\al)|^{\alpha}P(\al,t)]\nu_\alpha (dz).\tag {3.3}
\end{align}
Since $X\in (0,\infty)$, then $f_3(\al)=\epsilon\,\al\neq0$ and $|f_3(\al)|^\alpha \in C^2(R)$.

\begin{figure}[ht!]
\begin{subfigure}[b]{0.45\linewidth}
 \centering
  \includegraphics[width=\linewidth]{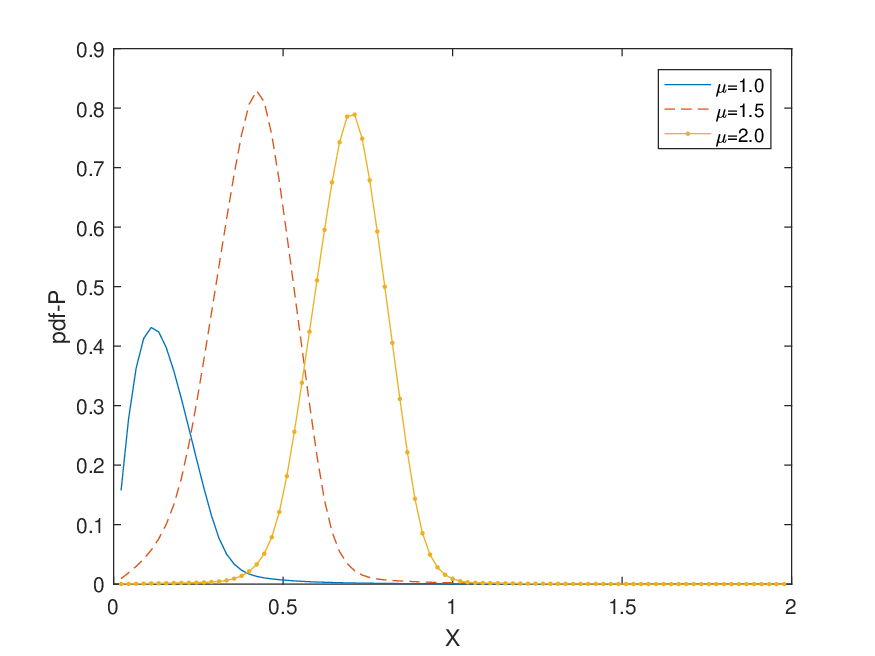}
  \centering
 \caption{Stationary probability density verse $\mu$.}
     \end{subfigure}
     \begin{subfigure}[b]{0.45\linewidth}
     \centering
  \includegraphics[width=\linewidth]{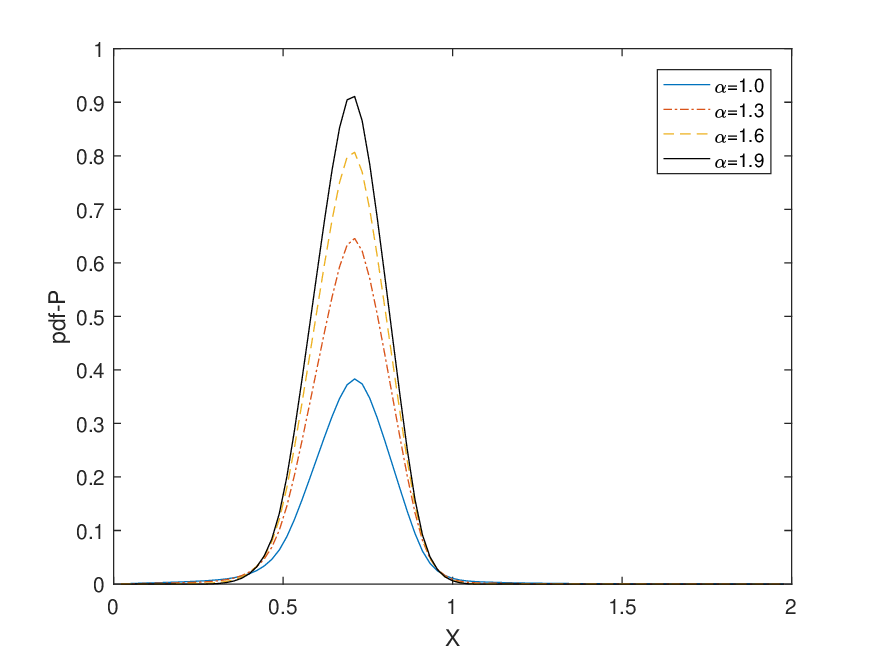}
 \caption{Stationary probability density verse $\alpha$.}
     \end{subfigure}
      \begin{subfigure}[b]{0.45\linewidth}
      \centering
  \includegraphics[width=\linewidth]{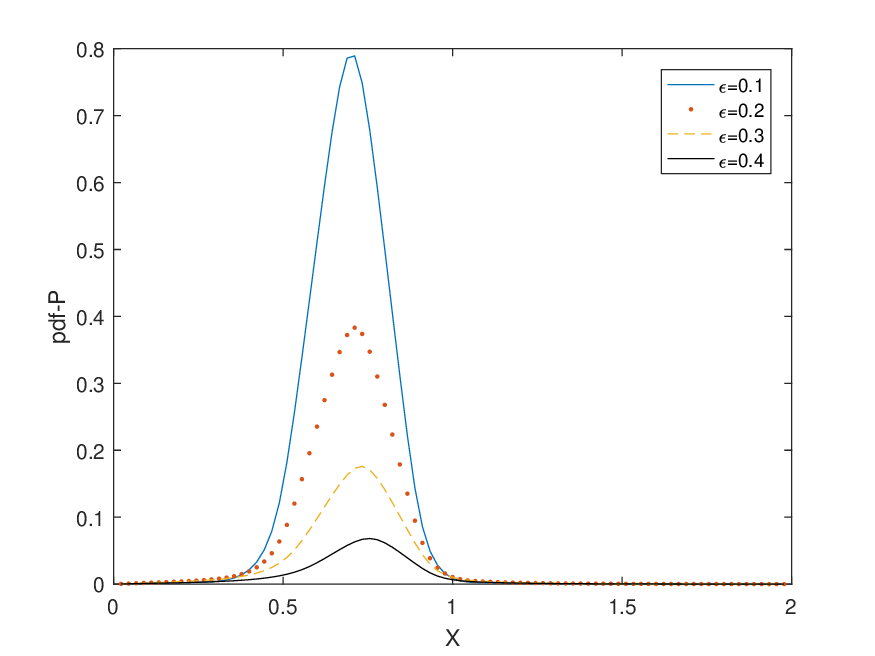}
 \caption{Stationary probability density verse $\epsilon$.}
     \end{subfigure}
\hfill
    \begin{subfigure}[b]{0.45\linewidth}
      \centering
  \includegraphics[width=\linewidth]{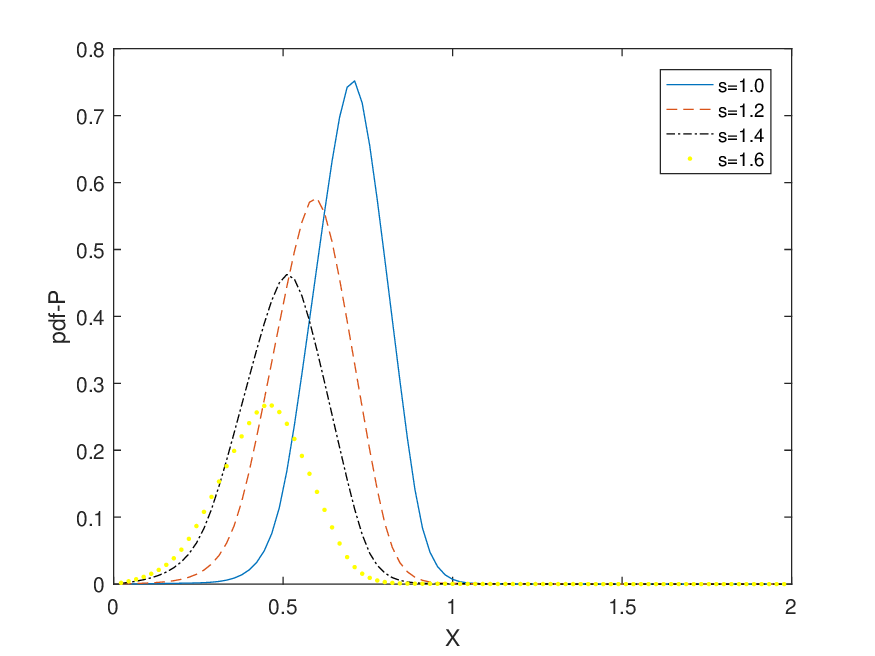}
 \caption{Stationary probability density verse $s$.}
     \end{subfigure}
     \caption{\label{FPE-M1} (Color online) Stationary probability density function for model (\ref{s1}). (a) $\alpha=1$, $\epsilon=0.1$, $s=1$. (b) $\mu=2$, $\epsilon=0.2$, $s=1$. (b) $\epsilon=1$, $\mu=2$, $s=1$. (c) $\alpha=1$, $s=1$, $\mu=2$. (d) $\epsilon=0.2$, $\alpha=1.5$, $\mu=2$. Parameters $\lambda=0$, $x_0=0.5$, time $T=50$.}

\end{figure}

To simulate the stationary density of the non-local Fokker-Planck equation in Eq. (\ref{p1}), we apply a numerical finite difference method developed in \cite{40}, and its numerical solutions as the parameters $\mu$, $\alpha$, $\epsilon$ and $s$ vary are plotted in Fig. \ref{FPE-M1}.

\section{Bifurcation analysis for a logistic growth model}\label{Sec4}
\subsection{\textbf{Deterministic saddle-node bifurcation}}\label{Sec41}
An ordinary differential equation (the noise free equation) describes the dynamics of Verhulst logistic equation \cite{41}:
 \begin{equation}\label{Det2}
\frac{dX}{dt}=s\,\al\, \left(1-\frac{\al}{M}\right), \quad \al_0=x_{0},\tag {4.1}
 \end{equation}
where the parameters $s$ and $M$, and the variable $X$ have the same biological meanings with the model (\ref{Det21}). Here the initial population size $x_0$ must be greater than zero.

\begin{figure}[ht!]
\centering
\begin{center}
  \includegraphics[width=0.5\linewidth]{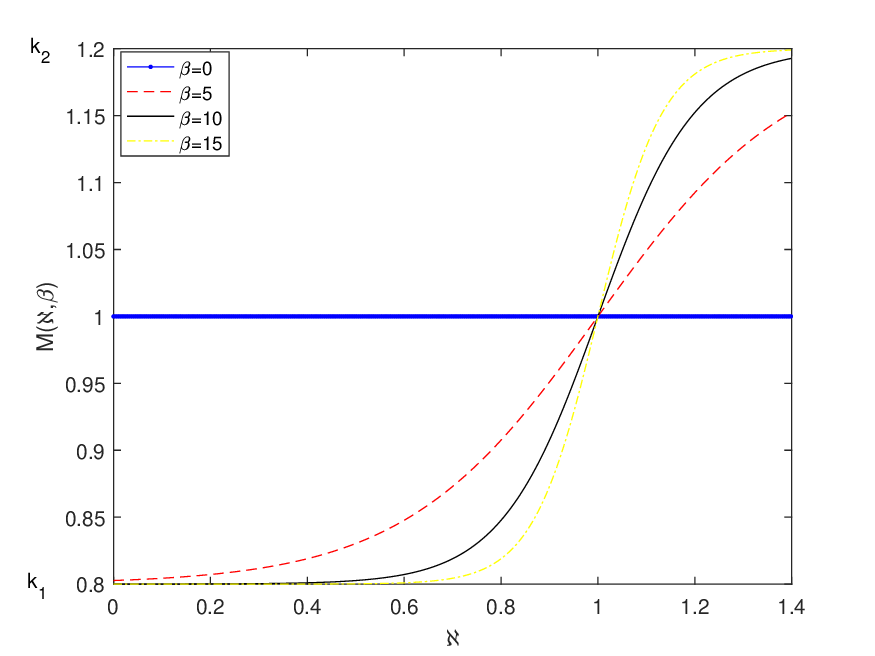}
 \caption{ (Color online) Phaseline diagram for carrying capacity in Eq. (\ref{carry}). Parameters $s=0.2$, $k_2=1.2$, $k_1=0.8$, and $\beta=0, 5, 10, 15$.}
 \label{carrying}
 \end{center}
\end{figure}

Instead of constant carrying capacity in Eq. (\ref{Det2}), we assume state-dependent carrying capacity,
\begin{equation}\label{carry}
M(\al)=\left(k_1+\frac{k_2-k_1}{1+e^{-\beta (\al-\varphi)}}\right),\tag {4.2}
\end{equation}
where $k_1,k_2 >0$, and $\beta >0$ controls the sensitivity of $M(\al)$ with respect to the change of the population size $\al$, and $\varphi$ is a constant given by $\varphi =\frac{k_1+k_2}{2}$. From Figure \ref{carrying}, we observe that $M(\al)$ lies between $k_1$ and $k_2$.
When $\beta$ is near bifurcation (discussed below), the effect of this form for $M(\al)$ is to replace the carrying capacity of a single point with an interval of carrying capacities where, in the absence of noise, $\dot{\al}(t)\sim 0$.

Replace $M$ in Eq. (\ref{Det2}) by $M(\al)$, yields,
 \begin{equation}\label{Det2M}
\frac{dX}{dt}=s\,\al\,\left(1-\frac{\al}{M(\al)}\right), \quad \al_0=x_{0}.\tag {4.3}
 \end{equation}
\begin{figure}[ht!]
\centering
      \begin{subfigure}[b]{0.5\textwidth}
        \centering
         \includegraphics[width=\textwidth]{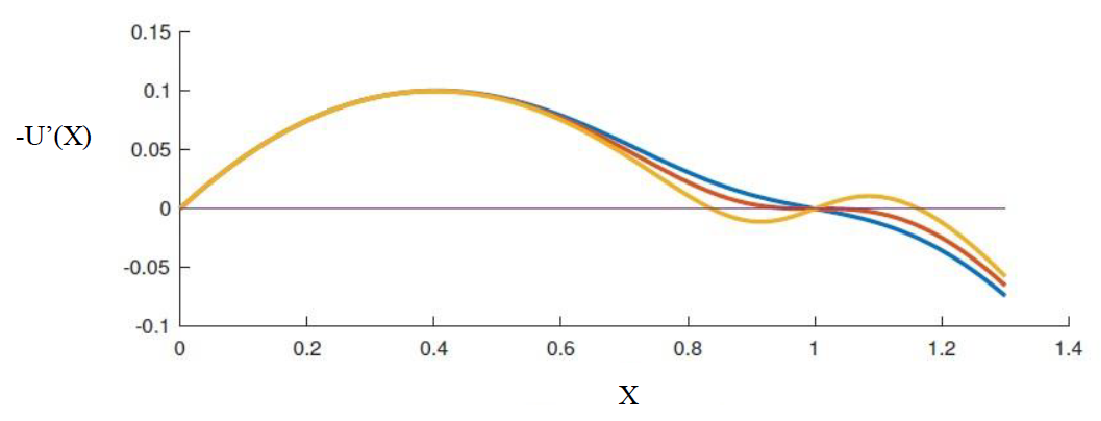}
         \caption{Phaseline diagram.}
          \end{subfigure}
     \begin{subfigure}[b]{0.5\textwidth}
        \centering
         \includegraphics[width=\textwidth]{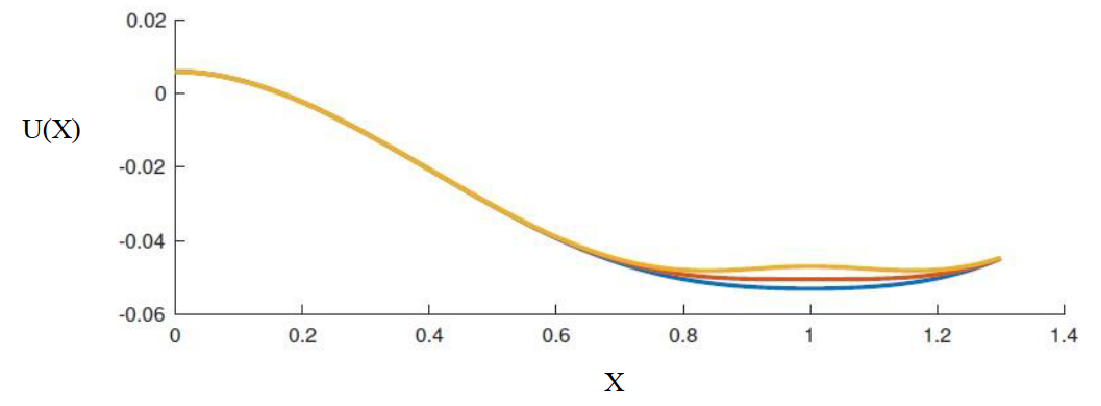}
         \caption{Potential function.}
         \end{subfigure}
         \caption{ \label{poten1}(Color online) Phaseline diagram and potential function $U$ for Eq. (\ref{Det2M}).}
\end{figure}
Equation (\ref{Det2M}) has two equilibrium solutions $\al_1(\beta) = 0$ and  $\al_2(\beta) = \varphi$, and it can be written as
  \begin{equation*}
  \frac{d\al}{dt}=-\frac{\partial U}{d\al},
  \end{equation*}
 where $U$ is the potential function defined by
  \begin{equation*}
  U(\al)=-\int s\,\al\left(1-\frac{\al}{M(\al)}\right)d\al.
  \end{equation*}
\begin{figure}[ht!]
    \centering
     \begin{subfigure}[b]{0.45\textwidth}
         \centering
         \includegraphics[width=\textwidth]{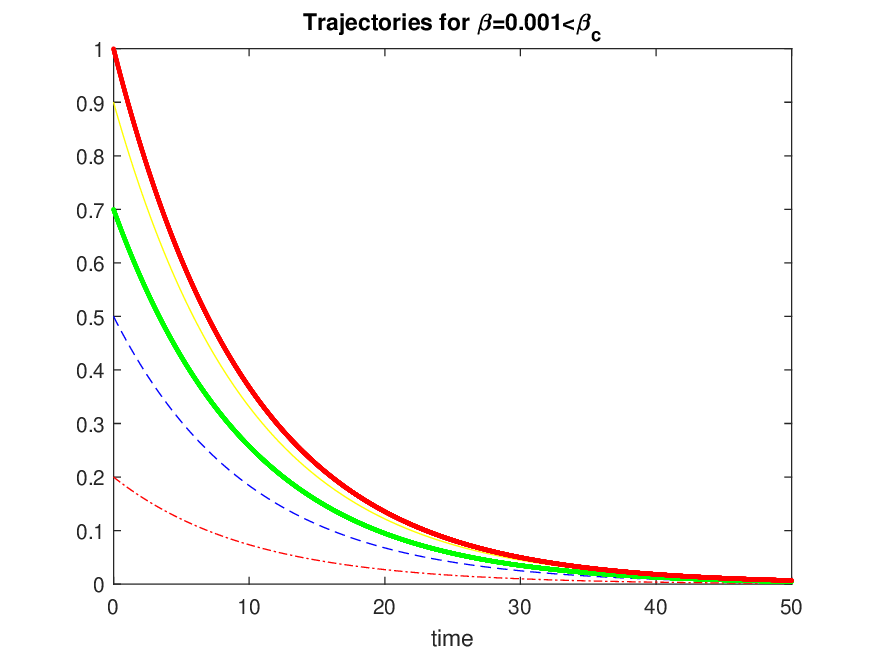}
         \caption{}
          \end{subfigure}
     \begin{subfigure}[b]{0.45\textwidth}
         \centering
         \includegraphics[width=\textwidth]{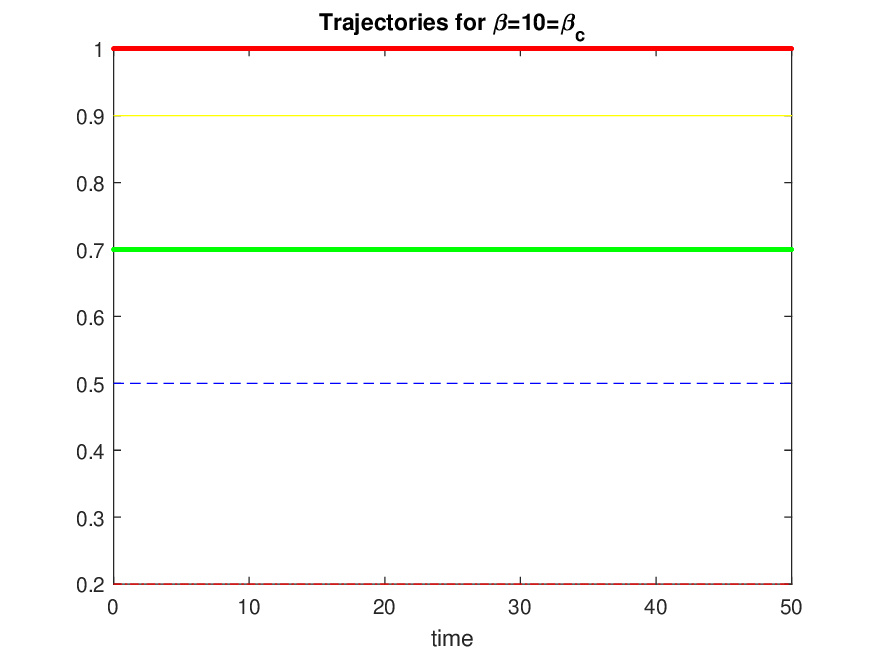}
         \caption{}
         \end{subfigure}
     \begin{subfigure}[b]{0.45\textwidth}
         \centering
         \includegraphics[width=\textwidth]{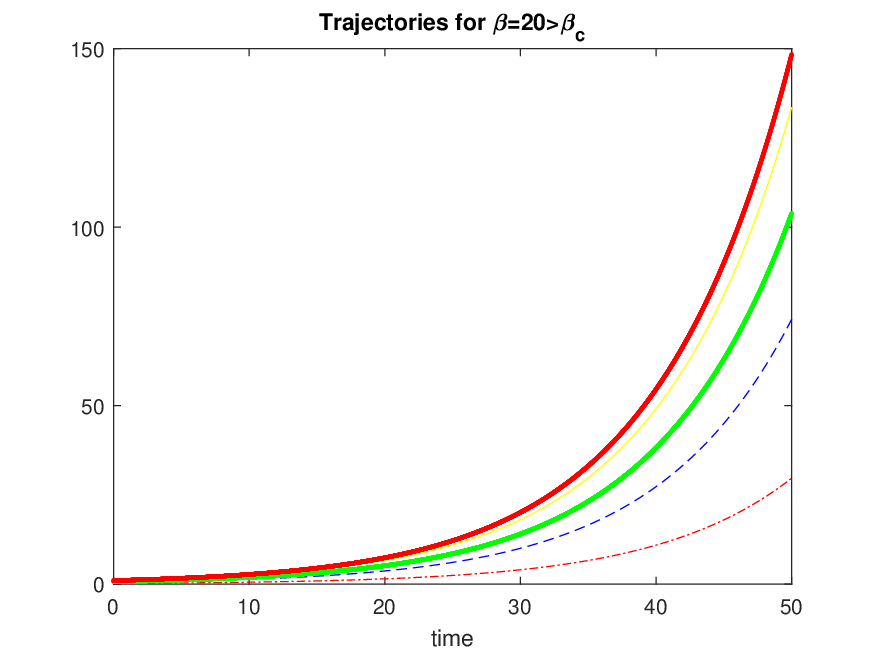}
         \caption{}
         \end{subfigure}
        \caption{\label{phase} (Color online) The solution of Eq.(\ref{sol}). (a) $\beta=0.001$. (b) $\beta=10$. (c) $\beta=20$. Parameters $s=0.2$, $k_2=1.2$, $ k_1=0.8$ and $u_0= 0.5, 0.7, 0.9, 1.0$.}

\end{figure}

The potential function $U$ exhibits two equilibrium states $\al_1(\beta) = 0$ and  $\al_2(\beta) = \varphi$. In biology, $\al_1$ tells the extinction state, where no population are present, and $\al_2$ indicates the state of stable population where population does not increase but stays at a certain constant.

Note that near a bifurcation point $\beta_C$ (discussed below), the potential wells shown in Fig. \ref{poten1}b have flat bottoms where $\frac{d\al}{dt}\sim 0$  in the absence of noise.

Substituting $\al=\varphi+u$ into Eq. (\ref{Det2M}) and dropping the nonlinear terms in the perturbation $u$ infer the linear approximation to Eq. (\ref{Det2M}) near $\al_2(\beta) = \varphi$,
\begin{equation}\label{delmod}
  \frac{du}{dt}=-s \left[1-\frac{1}{4} (k_2-k_1)\beta \right] u, \quad u(0)=u_0.\tag {4.4}
\end{equation}

The coefficient of $u$ on the right side of Eq. (\ref{delmod}) is negative when $\beta<\frac{4}{k_2-k_1}$ and positive when $\beta > \frac{4}{k_2-k_1}$, we conclude that
\begin{equation*}
  \beta_c=\frac{4}{k_2-k_1}
\end{equation*}
is a bifurcation value for Eq. (\ref{Det2M}). This bifurcation value is called saddle-node bifurcation.

Since $\beta_c=\frac{4}{k_2-k_1}$, we can rewrite Eq. (\ref{delmod}) as follows:
\begin{equation}\label{delmod1}
  \frac{du}{dt}=-s \left[1-\frac{\beta}{\beta_c} \right] u, \quad u(0)=u_0.\tag {4.5}
\end{equation}
The solution of Eq. (\ref{delmod1}) is given by
\begin{equation}\label{sol}
  u(t)=u_o\,e^{At},\tag {4.6}
\end{equation}
where $A=-s \left[1-\frac{\beta}{\beta_c} \right]$. Suppose that the initial value $u_0\neq 0$. Thus,

\[
\lim_{t\rightarrow\infty } u(t)=
\begin{cases}
\infty, \qquad  \qquad \beta_c<\beta,\\
u_0, \qquad \qquad \beta_c=\beta,\\
0, \qquad \qquad  \beta_c>\beta.\\
\end{cases}
\]
\subsection{\textbf{Stochastic bifurcation of a logistic growth model}}\label{Sec42}
\begin{figure}[ht!]
\begin{center}
  \includegraphics[width=0.45\linewidth]{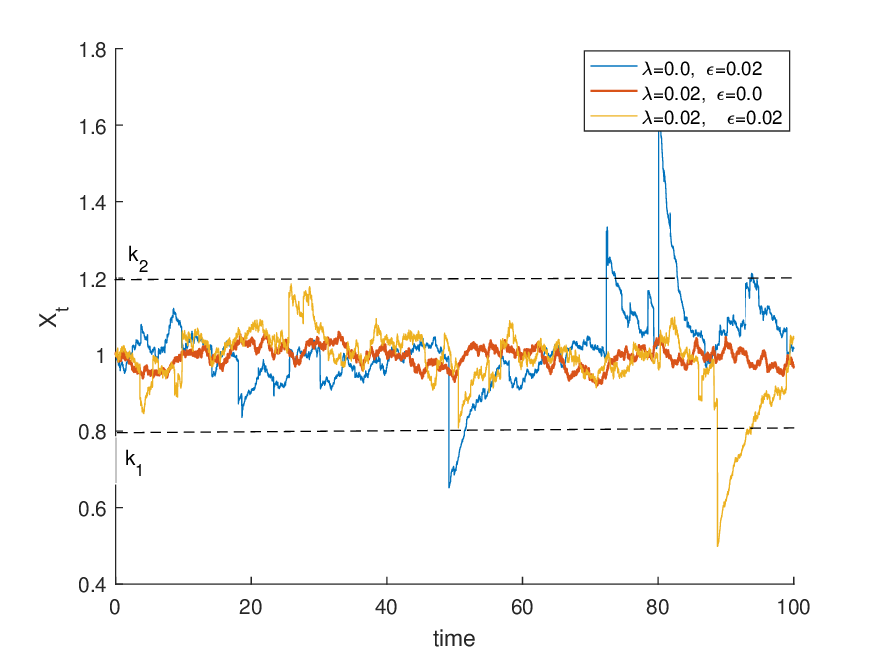}
 \caption{ Sample paths of Equation (\ref{SDEBMLM1}) for fixed $s=0.5,\, k_2=1.2,\,  k_1=0.8$, $x_0=1,$ $\alpha=1.5$. and $\beta=5$.}
 \label{sample-M2}
 \end{center}
\end{figure}

In this subsection, we consider the following stochastic model:
\begin{equation}\label{SDEBMLM1}
  d\al =s\,\al\left(1-\frac{\al}{M(\al)}\right)dt+\lambda \al\,dB_t+\epsilon \al\,dL_{t}^{\alpha},\quad \al_0=x_0,\tag {4.7}
  \end{equation}
where $B_t$ is the standard Brownian motion, $\lambda$ and $\epsilon$ represent the intensities of Gaussian noise and non-Gaussian noise respectively, $L_{t}^{\alpha}$ is an $\alpha$-stable L\'evy motion. Here the noises are multiplicative because the diffusion and the intensity coefficient in (\ref{SDEBMLM1}) depend on $\al$. Figure \ref{sample-M2} shows a sample path of the solution of the SDE in (\ref{SDEBMLM1}) when $\epsilon=0$.\\
The generator $\mathcal{A}$ for the model (\ref{SDEBMLM1}), or for the solution process $\al$, is
\begin{equation*}
\mathcal{A}\varphi(\al)=\left[s\,\al\,\left(1-\frac{\al}{M(\al)}\right)\right]\varphi'(\al)+\frac{1}{2}\lambda^{2}\al^{2}\varphi''(\al)+\int_{\mathbb{R}\setminus\{0\}}\big(\varphi(\al+\epsilon\,\al y)-\varphi(\al)\big)\nu_{\alpha}(dy),
\end{equation*}
where $\varphi$ belongs to the domain of definition for operator $\mathcal{A}$.

Setting $z= \epsilon\,\al\,y$, the probability density $P(\al,t)$ of the Fokker-Planck equation for the model (\ref{SDEBMLM1}) for the solution process $\al$ with initial condition $x_0$  and $P(\al,0)=\sqrt {\frac{40}{\pi}}e^{-40x_0^2}$, becomes the following non-local partial differential equation:
\begin{align}\label{fpe-F2}
\partial_{t}P(\al,t)&=- \partial_\al\left[s\,\al\,\left(1-\frac{\al}{M(\al)}\right)P(\al,t)\right]+\frac{1}{2}\partial_{\al\al}\left[\lambda^{2}\al^{2}P(\al,t)\right]\nonumber\\
&+|\epsilon\,\al|^{\alpha}\,\int_{\mathbb{R}\setminus\{0\}}[|f_3(\al+z)|^{\alpha}P(\al+z,t)-|f_3(\al)|^{\alpha}P(\al,t)]\nu_\alpha(dz).\tag {4.8}
\end{align}
To find numerically solution of the probability density $P(\al,t)$ in the non-local Fokker-Planck equation (\ref{fpe-F2}), we apply a numerical finite difference method used in  \cite{40}.
 \begin{figure}[ht!]
 \centering
\begin{subfigure}[b]{0.45\linewidth}
\centering
  \includegraphics[width=\linewidth]{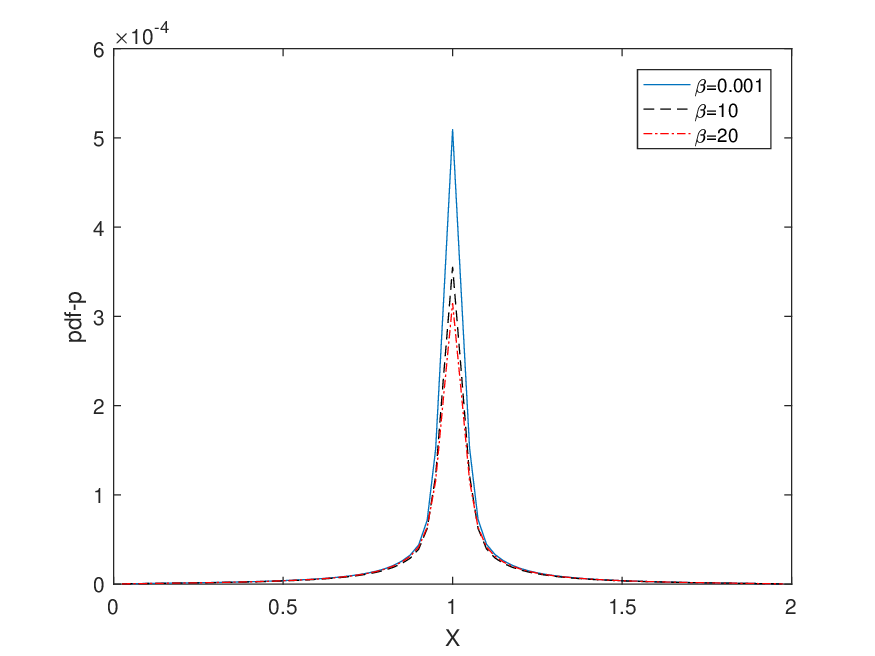}
 \caption{Stationary probability density verse $\beta$.}
     \end{subfigure}
     \begin{subfigure}[b]{0.45\linewidth}
     \centering
  \includegraphics[width=\linewidth]{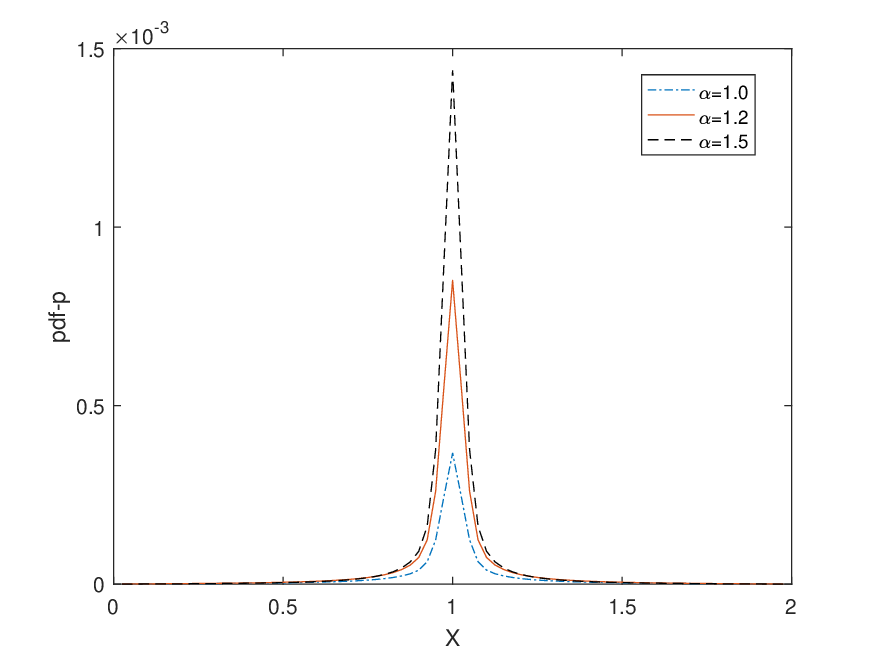}
 \caption{Stationary probability density verse $\alpha$.}
     \end{subfigure}
     \hfill
      \begin{subfigure}[b]{0.45\linewidth}
      \centering
  \includegraphics[width=\linewidth]{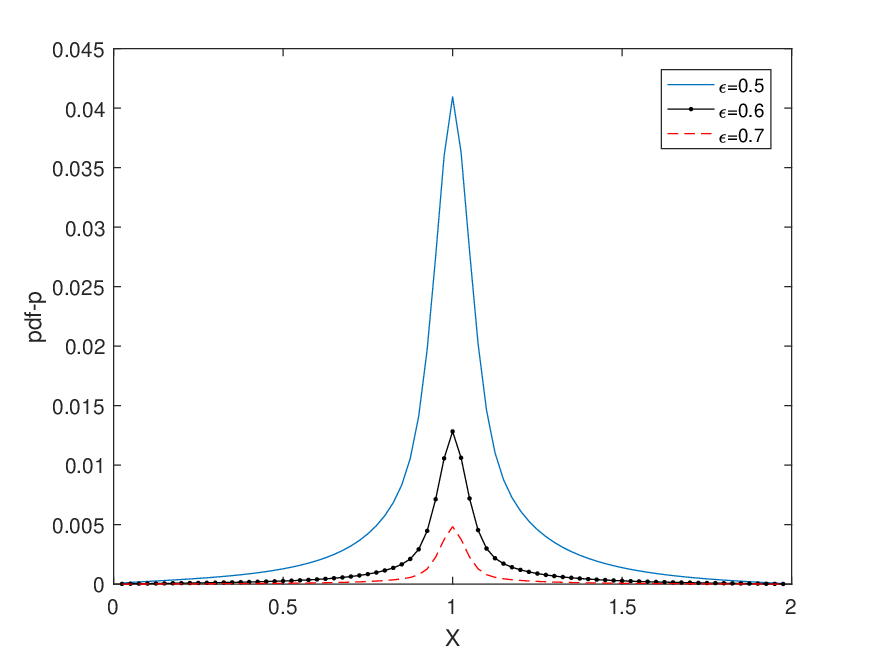}
 \caption{Stationary probability density verse $\epsilon$.}
     \end{subfigure}
      \begin{subfigure}[b]{0.45\linewidth}
      \centering
  \includegraphics[width=\linewidth]{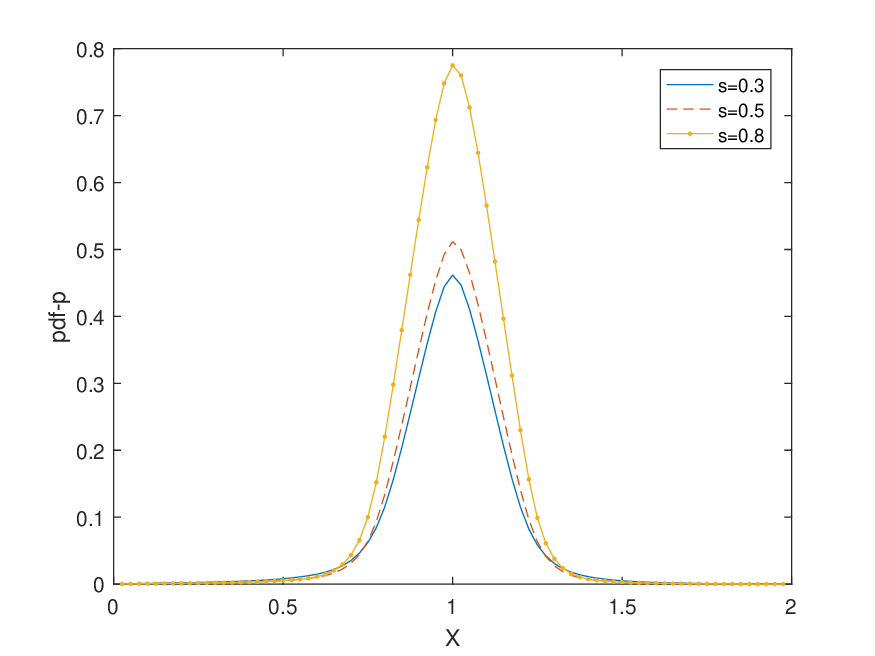}
 \caption{Stationary probability density for varying growth
rates $s$.}
     \end{subfigure}
 \caption{ \label{FPE-M2}(Color online) Stationary probability density function for model (\ref{SDEBMLM1}). (a) $\alpha=1$, $\epsilon=1$, $s=0.1$. (b) $\beta=8$, $\epsilon=1$,  $s=0.1$. (c) $\alpha=1$, $\beta=8$, $s=0.1$. (d) $\alpha=1$, $\epsilon=0.3$, $\beta=8$. Parameters $\lambda=0$, $x_0=0.5$, time $T=50$.}
\end{figure}
\section{Numerical experiments}\label{Sec5}
This section describes our numerical findings of the above deterministic models (\ref{deter}) and (\ref{Det2M}), and  stochastic models given in Eqs. (\ref{s1}) and (\ref{SDEBMLM1}) one by one. In order to get the critical parametric behaviors for stochastic P-bifurcation, we discuss the parametric impacts on stochastic P-bifurcation of the models by finding the peaks in the stationary density function curve.
\subsection{{Results for a growth model}}
The phase line diagram of the deterministic model (\ref{deter}) is plotted in Fig. \ref{Ph-M1}a for various of $\mu$. From this, we can see that the model has two equilibria points when $\mu>1$. While for $\mu \leq 1$ it has only one. Fig. \ref{Ph-M1}b shows the bifurcation diagram of Eq. (\ref{deter}). If $\mu<1$ there are two fixed points: $\al_1=0$ is stable ( solid black line) and  $\al_2=\frac{1}{s}\ln\mu$ is unstable (broken green line). The two points coincide at $(1,0)$. If $\mu>1$, the stability changes $\al_1$ is unstable ( broken black line) and $\al_2$ is stable (solid green line).  However, the stability of the fixed point along a given curve changes on passing through bifurcation value $\mu=1$. Thus this type of bifurcation is known as transcritical bifurcation.

Figure \ref{sample-M1} shows the numerical solution ( sample paths) of model (\ref{s1}) with $q-p-\frac{\lambda^2}{2}>0$.

The graphs of the probability density function $P$ for the Fokker-Planck equation for the solution of the model (\ref{s1}) are presented in  Figure \ref{FPE-M1} for different values of the parameter $\mu$, $\alpha$, $\epsilon$ and the growth rate $s$ of the system. Let us discuss the effect of these parameters in the model.\\
\textbf{\emph{Impact of the parameter $\mu$}}: Figure \ref{FPE-M1}a shows the probability density function $P$ with respect to $\mu$. Here $P$ evolves from the unimodal to the flatter kurtosis shape, and the maximum of the density functions move from the center (free-state of the population) to the right side of the domain ( stable state of the population).\\
\textbf{\emph{Effect of the parameter $\alpha$}}: The probability density function verse $\alpha$ is plotted in Figure \ref{FPE-M1}b. This shows that as $\alpha$ value increases, the peak of the probability density function also increases.\\
\textbf{{\emph{Influence of the parameter $\epsilon$}}}: The shape of the probability density function changes from unimodal to flatter as $\epsilon$ value increases. We can observe this phenomenons from Figure \ref{FPE-M1}c. \\
\textbf{\emph{Effect of the parameter $s$}}: In Figure \ref{FPE-M1}d we plot the probability density function with respect to the per capita daily adult mortality rate $s$.

Numerical simulations in Figure \ref{FPE-M1} imply that the system (\ref{s1}) is sensitive to small changes in the parameter values. Furthermore, small values of $\epsilon$, $s$ and $\mu$  indicate a transition in the L\'evy noise case.
\subsection{{Results for the logistic growth model}}
In Fig. \ref{carrying}a we present the phase line diagram for the state-dependent carrying capacity in Eq. (\ref{carry}) with different values of $\beta$. As we can see from this Figure that all curves are intersected at the point $(1,1)$.

Figure \ref{phase} shows the trajectories of Eq. (\ref{sol}) at various values of $\beta$ and initial value $u_0$. In Fig. \ref{phase}a when the value of $\beta$ is less than $\beta_C$, the solution of Eq. (\ref{sol}) goes to zero. Fig. \ref{phase}b shows if $\beta=\beta_C$, $u(t)=u_0$. While for $\beta > \beta_C$, the solution $u(t)$ goes to $\infty$.

The potential function and phase line diagram of the deterministic model (\ref{Det2M}) are presented in Figure \ref{poten1}. As can be seen from Figure \ref{poten1}b, the potential functions have flat-bottom near the bifurcation when $\frac{d\al}{dt}\sim 0$, and it shows mono-stable potential well. The potential function has one maximum corresponds to the unstable state and one minimum corresponds to the stable state of the population.

Figure \ref{sample-M2} plots the result of a numerical simulation of the model (\ref{SDEBMLM1}). From this Figure we can see that the sample paths are positive solutions and also bounded below by $k_1$ and bounded above by $k_2$. Here parameters $s=0.5$, $k_1=0.8$, $k_2=1.2$, $\alpha=1.5,$ and $\beta=5$ are fixed.

The numerical solution of the probability density function $P$ for the FPE of the solution of the system in (\ref{SDEBMLM1}) is presented in Fig. \ref{FPE-M2} with fixed values of $\lambda=0$ and $x_0=0.5.$\\
\textbf{\emph{Impact of the parameter $\beta$}}: In Fig. \ref{FPE-M2}a, we plot the probability density function $P$ at different values of $\beta$. From this Figure, we observe that as $\beta$ increases, $P$ decreases.\\
\textbf{\emph{Impact of the parameter $\alpha$}}: The stationary probability density functions of the system (\ref{SDEBMLM1}) are plotted in Fig. \ref{FPE-M2}b for various values of $\alpha$. For larger stability index, the larger stationary probability density. \\
\textbf{\emph{Dependence of different $\epsilon$}}:
 Next, we observe the dependence of the density function $P$ for different values of noise intensity $\epsilon$ which are plotted in Fig. \ref{FPE-M2}c. The larger $\epsilon$ produces small $P$, i.e., as the noise intensity $\epsilon$ increases, the stationary probability density decreases and its existing peak decreases in height.\\
\textbf{\emph{Impact of the parameter $s$}}: From Figure \ref{FPE-M2}d, we can see that the probability density function $P$ evolves from higher kurtosis to the lower one, as $s$ value increases. It means that for a larger value of $s$, the shape of $P$ becomes flatter.

Comparing with the results in Figure \ref{FPE-M2}, we can observe that it needs a small value of noise intensity $\epsilon$, $s$ and $\beta$ to induce a transition in the non-Gaussian case.
\section{Conclusion and discussion}\label{Sec6}
Bifurcation analysis is a key tool for the analysis of dynamic systems and it also enables the researchers to qualitatively study the behavior of trajectories without utilizing the analytical solution of the underlying ordinary differential equations (ODEs) and stochastic differential equations (SDEs).

In the present study, we investigated the stochastic P-bifurcations of two stochastic population models (i) a growth model with state dependent growth rate and constant death rate, and (ii) a logistic growth model with state dependent carrying capacity, both of which are driven by multiplicative symmetric $\alpha$-stable L\'evy noise.  We provided a new approach to study the behavior of the above two models. We reviewed the deterministic bifurcation of the models. The effects of parameters $\mu$, $s$, $\beta$, $\alpha$, and $\epsilon$ on the population models are analyzed. Since the stochastic bifurcations are about to the qualitative changes of the stationary probability density $P$ of the Fokker-Planck equation which is a non-local partial differential equation on the stochastic system with an $\alpha$-stable L\'evy motion. Here, we used a MATLAB program to generate the stochastic P-bifurcations of our two systems (\ref{s1}) and (\ref{SDEBMLM1}).

In growth model (\ref{deter}) when the birth rate $q$ is larger than the death rate $p$ ($\mu>1$) there are two fixed points;  one of the fixed points is stable (represented by the solid green color), and the other fixed point is unstable (represented by the broken black color). While the birth rate $q$ is smaller than the death rate $p$, the stability changed. The two fixed points $\al_1$ and $\al_2$ coincide when the birth rate $q$ is equal to the death rate $p$ ( $\mu=1$). In the biological point of view, the birth rate $q$ must be greater than the death rate $p$, so at $\al=0$ extinction can occur.

In logistic model (\ref{Det21-SDE}), we assumed state-dependent carrying capacity. The deterministic counterparts of the systems (\ref{s1}) and (\ref{SDEBMLM1}) have transcritical bifurcation and saddle-node bifurcation, respectively. And from our numerical experiment, we observed that the change in the values of the parameters causes a change of the shape of the stationary probability density. We also analyzed the parametric influences on stochastic P-bifurcation of the system. The figures are plotted in Fig. \ref{FPE-M1} and Fig. \ref{FPE-M2}.

A peak appeared at the middle of the domain, which means that a transition occurs. From a biological perspective, the extinction of the population can not occur. From Fig. \ref{FPE-M2}, we observed that the maximum of the stationary density function remains at the center of the domain. This showed the probability of the  population growth steady at this point becomes smaller. This does not imply a decrease in the total population size. A slight change in the value of the parameter has a crucial effect in the behavior of the population systems.
\section*{Acknowledgments}
The authors are grateful to Professor Bernardo Spagnolo and Professor Alexander A. Dubkov
for thoughtful comments and suggestions on stochastic bifurcations of logistic growth models driven by L\'evy motions. The authors acknowledge support from the NSFC grant 12001213.

\end{document}